\documentclass[12pt,a4paper]{article}
\usepackage{amsmath}
\usepackage{amssymb}
\usepackage[thmmarks,amsmath]{ntheorem}

\makeatletter
    
    \newcommand{\Rmnum}[1]{\expandafter\@slowromancap\romannumeral #1@}
  \makeatother

\newtheorem{propo}{Proposition}[section]

\newtheorem{lemma}{Lemma}[section]

\newtheorem{theo}[propo]{Theorem}

\numberwithin{equation}{section}
\begin{document}
\title{An Upper Bound for the Number of Solutions of Ternary Purely Exponential Diophantine Equations }

\author{ Yongzhong Hu and Maohua Le }
\maketitle
\maketitle \edef \tmp {\the \catcode`@}
   \catcode`@=11
   \def \@thefnmark {}
   \@footnotetext {Supported by the National Natural Science Foundation of China(No.10971184)}
   \catcode`@=\tmp
   \let\tmp = \undefined

\begin{abstract}   Let $a,b,c$ be fixed coprime positive integers with $\min\{a,b,c\}>1$. In this paper, combining the Gel'fond-Baker method with an elementary approach, we prove that if $\max\{a,b,c\}>5\times 10^{27}$, then the equation $a^x+b^y=c^z$ has at most three positive integer solutions $(x,y,z)$.
\end{abstract}



\section {Introduction}
\quad  Let $\mathbb{N}$ be the set of all positive integers. Let $a,b,c$ be fixed coprime positive integers with $\min\{a,b,c\}>1$.
In 1933, K. Mahler \cite{mah14} used his $p-$adic analogue of the Thue-Siegel  method to prove that the ternary purely exponential diophantine equation
\begin{equation}\label{1.1}
a^x+b^y=c^z, x,y,z\in\mathbb{N}
\end{equation}
has only finitely many solutions $(x,y,z)$. His method is ineffective. An effective result for solutions of $(\ref{1.1})$ was given by A.O. Gel'fond \cite{gel7}. Let $N(a,b,c)$ denote the number of solutions $(x,y,z)$ of $(\ref{1.1})$. As a straightforward consequence of an upper bound for the number of solutions of binary $S-$unit equations due to F. Beukers and H. P. Schlickewei \cite{beu2}, we have $N(a,b,c)\leq2^{36}$. In recent years, many papers investigated the exact values of $N(a,b,c)$. The known results showed that $(\ref{1.1})$ has only a few solutions for some special cases(see \cite{cip5}, \cite{fuj6}, \cite{lem12}, \cite{luc13}, \cite{miy15}, \cite{miy16}, \cite{miy17}, \cite{miy18}, \cite{miy19}, \cite{miy20}, \cite{miy21} and \cite{ter22}). Very recently, the authors \cite{hu8} proved that if $a,b,c$ satisfy certain divisibility conditions and $\max\{a,b,c\}$ is large enough, then $(\ref{1.1})$ has at most one solution $(x,y,z)$ with $\min\{x,y,z\}>1$. In this paper we prove a general result as follows:
\begin{theo}
If $\max\{a,b,c\}>5\times10^{27}$, then $N(a,b,c)\leq3$.
\end{theo}

Notice that if $(a,b,c)=(3,5,2)$, then $(\ref{1.1})$ has exactly three solutions $(x,y,z)=(1,1,3),(3,1,5)$ and $(1,3,7)$. Perhaps, in general, $N(a,b,c)\leq3$ is the best upper bound for $N(a,b,c)$.

\section {An upper bound for the solutions of $(\ref{1.1})$}

\quad  In \cite{hu8}, combining a lower bound for linear forms in two logarithms and an upper bound for the $p-$adic logarithms due to M. Laurent \cite{lau9} and Y. Bugeaud \cite{bug3} respectively, the authors proved that all solutions $(x,y,z)$ of $(\ref{1.1})$ satisfy $\max\{x,y,z\}<155000(\log\max\{a,b,c\})^3$, where $\log$ is used for natural logarithm. In this section, using the same method as in \cite{hu8}, we make a slight improvement as follows:

\begin{theo}
All solutions $(x,y,z)$ of $(\ref{1.1})$ satisfy
\begin{equation}\label{2.1}
\max\{x,y,z\}<6500(\log\max\{a,b,c\})^3.
\end{equation}
\end{theo}

The proof of Theorem 2.1 depends on the following lemmas.


\begin{lemma}\label{2l2}{\rm(\cite{lau10}, Corollaire 2 et Tableau 2)}.
 Let $\alpha_1,\alpha_2,\beta_1,\beta_2$ be positive integers with $\min\{\alpha_1,\alpha_2\}\geq2$. Further let $\Lambda=\beta_1\log\alpha_1-\beta_2\log\alpha_2$. If $\Lambda\not=0$, then
 $$\log|\Lambda|>-32.31(\log\alpha_1)(\log\alpha_2)(\max\{10,0.18+\log(\frac{\beta_1}{\log\alpha_2}+\frac{\beta_2}{\log\alpha_1})\})^2.$$
\end{lemma}

\begin{lemma}.\label{2l3}
 Let $\alpha_1,\alpha_2$ be odd integers with $\min\{|\alpha_1|,|\alpha_2|\}\geq3$, and let $\beta_1,\beta_2$ be positive integers. Further let $\Lambda^\prime=\alpha_1^{\beta_1}-\alpha_2^{\beta_2}$. If $\Lambda^\prime\not=0$ and $\alpha_1\equiv\alpha_2\equiv1\pmod4$, then
 $${\rm ord}_2\Lambda^\prime<19.57(\log|\alpha_1|)(\log|\alpha_2|)(\max\{12\log2,$$
 $$0.4+\log(2\log2)+\log(\frac{\beta_1}{\log|\alpha_2|}+\frac{\beta_2}{\log|\alpha_1|})\})^2,$$
 where ${\rm ord_2}\Lambda^\prime$ is the order of $2$ in $|\Lambda^\prime|$.
\end{lemma}


{\bf {\it Proof.}}\,\,
This is the special case of [3, Theorem 2] for $p=2, y_1=y_2=1,\alpha_1\equiv\alpha_2\equiv1\pmod4,g=1$ and $E=2$.$\hfill{} \Box$\\

Throughout this section, let $(x,y,z)$ be a solution of $(\ref{1.1})$. Since $\max\{a^x,b^y\}<c^z$, we have

\begin{equation}\label{2.2}
\max\{x\log a,y\log b\}<z\log c.
\end{equation}

\begin{lemma}\label{2l4}.
If $\min\{a^{2x},b^{2y}\}<c^z$, then

\begin{equation}\label{2.3}
\max\{x,y,z\}<4663(\log\max\{a,b,c\})^2.
\end{equation}
\end{lemma}
{\bf {\it Proof.}}\,\,By the symmetry of $a^x$ and $b^y$ in $(\ref{1.1})$, it suffices to prove the lemma for the case that $\min\{a^{2x},b^{2y}\}=a^{2x}<c^z$. Then we have

\begin{equation}\label{2.4}
2x\log a<z\log c.
\end{equation}
Since $a^x<b^y$ and $c^z\geq a+b\geq5$, we have $a^x/(2b^y+a^x)=a^x/(b^y+c^z)<2a^x/3c^z<2/3c^{z/2}\leq2/3\sqrt{5}$ and
$$z\log c=\log(b^y+a^x)=y\log b+\log(1+\frac{a^x}{b^y})$$

\begin{equation}\label{2.5}
<y\log b+\frac{a^x}{b^y}<y\log b+\frac{2a^x}{a^x+b^y}<y\log b+\frac{2}{c^{z/2}}.
\end{equation}

Let $(\alpha_1,\alpha_2,\beta_1,\beta_2)=(c,b,z,y)$ and $\Lambda=\beta_1\log\alpha_1-\beta_2\log\alpha_2$. By $(\ref{2.5})$, we get
$0<\Lambda<2/c^{z/2}$ and

\begin{equation}\label{2.6}
\log2-\log\Lambda>\frac{z}{2}\log c.
\end{equation}
Since $\min\{\alpha_1,\alpha_2\}\geq2$, using Lemma $\ref{2l2}$, we have

\begin{equation}\label{2.7}
\log\Lambda>-32.31(\log c)(\log b)(\max\{10,0.18+\log(\frac{z}{\log b}+\frac{y}{\log c})\})^2.
\end{equation}

When $10\geq0.18+\log(z/\log b+y/\log c)$, by $(\ref{2.6})$ and $(\ref{2.7})$, we get $\log2 +3231(\log c)(\log b)>z(\log c)/2$ and

\begin{equation}\label{2.8}
z<6464\log b.
\end{equation}
Further, by $(\ref{2.4})$ and $(\ref{2.8})$, we obtain

\begin{equation}\label{2.9}
y<6464\log c,x<\frac{6464(\log b)(\log c)}{2\log a}<4663(\log b)(\log c),
\end{equation}
since $a\geq2$. Furthermore, since $\max\{a,b,c\}\geq5$, we see from $(\ref{2.8})$ and $(\ref{2.9})$ that $(\ref{2.3})$ holds.

When $10<0.18+\log(z/\log b+y/\log c)$, by $(\ref{2.2})$, $(\ref{2.6})$ and $(\ref{2.7})$,

\begin{equation}\label{2.10}
\log2+32.31(\log c)(\log b)(0.18+\log(\frac{2z}{\log b}))^2>\frac{z}{2}\log c,
\end{equation}
whence we obtain
$$\frac{z}{\log b}<\frac{2\log2}{(\log b)(\log c)}+64.62(0.88+\log(\frac{z}{\log b}))^2$$
\begin{equation}\label{2.11}
<2+64.62(0.88+\log(\frac{z}{\log b}))^2.
\end{equation}
Let $t=z/\log b$ and $F(t)=t-64.62(0.88+\log t)^2-2$. Then we have $F(6000)>0$ and $F^\prime(t)=1-129.24(0.88+\log t)/t>0$ for $t\geq6000$. Hence, we get $F(t)>0$ for $t\geq 6000$. Therefore, we see from $(\ref{2.11})$ that $z$ satisfies $(\ref{2.8})$, and hence, $x$ and $y$ satisfy $(\ref{2.9})$. Thus, the lemma is proved.$\hfill{} \Box$\\

{\bf {\it Proof of Theorem 2.1.}}\,\,By Lemma $\ref{2l4}$, the theorem holds if $\min\{a^{2x},b^{2y}\}<c^z$. We may therefore assume that $\min\{a^{2x},b^{2y}\}>c^z$. Then we have
\begin{equation}\label{2.12}
z\log c<\min\{2x\log a,2y\log b\}.
\end{equation}
Further, by $(\ref{2.2})$ and $(\ref{2.12})$, the theorem holds if $\min\{x,y,z\}=1$. We may assume that $\min\{x,y,z\}>1$.

We first consider the case that $2|a$. Since $x>1$, by $(\ref{1.1})$, we have $2\not|bc,$$c^z-b^y\equiv a^x\equiv0\pmod4,(-1)^{(c-1)z/2}\equiv c^z\equiv b^y\equiv(-1)^{(b-1)y/2}\pmod4 $ and

\begin{equation}\label{2.13}
(-1)^{(c-1)z/2}=(-1)^{(b-1)y/2}.
\end{equation}
Let $(\alpha_1,\alpha_2,\beta_1,\beta_2)=((-1)^{(c-1)/2}c,(-1)^{(b-1)/2}b,z,y)$ and $\Lambda^\prime=\alpha_1^{\beta_1}-\alpha_2^{\beta_2}$. Then, $\alpha_1$ and $\alpha_2$ are odd integers satisfy
\begin{equation}\label{2.14}
\min\{|\alpha_1|,|\alpha_2|\}\geq3,\alpha_1\equiv\alpha_2\equiv1\pmod4.
\end{equation}
By  $(\ref{2.13})$, we have $\Lambda^\prime=(-1)^{(c-1)z/2}(c^z-b^y)=(-1)^{(c-1)z/2}a^x$, whence we get $\Lambda^\prime\not=0$ and

\begin{equation}\label{2.15}
{\rm ord_2}\Lambda^\prime\geq x.
\end{equation}
On the other hand, using Lemma $\ref{2l3}$, we have
$${\rm ord_2}\Lambda^\prime<19.57(\log c)(\log b)(\max\{12\log2,$$
\begin{equation}\label{2.16}
0.4+\log(2\log2)+\log(\frac{z}{\log b}+\frac{y}{\log c})\})^2.
\end{equation}

When $12\log2\geq0.4+\log(2\log2)+\log(z/\log b+y/\log c)$, we have

\begin{equation}\label{2.17}
z<e^{12\log2-0.4-\log(2\log2)}\log b<2000\log b.
\end{equation}
Therefore, by  $(\ref{2.2})$ and $(\ref{2.17})$, $(\ref{2.1})$ holds.

When $12\log2<0.4+\log(2\log2)+\log(z/\log b+y/\log c)$, by $(\ref{2.15})$ and $(\ref{2.16})$, we have

\begin{equation}\label{2.18}
x<19.57(\log c)(\log b)(0.74+\log(\frac{z}{\log b}+\frac{y}{\log c}))^2.
\end{equation}
Further, by  $(\ref{2.2})$, $(\ref{2.12})$ and $(\ref{2.18})$, we get
\begin{equation}\label{2.19}
\frac{z\log c}{2\log a}<x<19.57(\log c)(\log b)(0.74+\log(\frac{2z}{\log b}))^2,
\end{equation}
whence we obtain
\begin{equation}\label{2.20}
\frac{z}{\log b}<39.14(\log a)(1.44+\log(\frac{z}{\log b}))^2.
\end{equation}

Let $t=z/\log b$ and $F(t)=t-39.14(\log a)(1.44+\log t)^2$. Since $a\geq2$,we have $F(6500(\log a)^2)>0$ and $F^\prime(t)=1-78.28(\log a)(1.44+\log t)/t>0$ for $t\geq6500(\log a)^2$. Hence, we get $F(t)>0$ for $t\geq 6500(\log a)^2$. Therefore, we see from $(\ref{2.20})$ that

\begin{equation}\label{2.21}
z<6500(\log a)^2(\log b).
\end{equation}
Further, by $(\ref{2.2})$ and $(\ref{2.21})$, we have

\begin{equation}\label{2.22}
x<6500(\log a)(\log b)(\log c),y<6500(\log a)^2(\log c).
\end{equation}
By  $(\ref{2.21})$ and $(\ref{2.22})$, $\max\{x,y,z\}$ satisfies $(\ref{2.1})$ and the theorem is true if $2|a$.

By the symmetry of $a$ and $b$ in  $(\ref{1.1})$, using the same method as in the above proof, we can prove that the theorem is true if $2|b$.

 We finally consider the case that $2|c$. Since $z>1$, we have $a^x+b^y\equiv c^z\equiv0\pmod4, a^x\equiv-b^y\pmod4$ and

 \begin{equation}\label{2.23}
(-1)^{(a-1)x/2}=-(-1)^{(b-1)y/2}.
\end{equation}
 Let $(\alpha_1,\alpha_2,\beta_1,\beta_2)=((-1)^{(a-1)/2}a,(-1)^{(b-1)/2}b,x,y)$ and $\Lambda^\prime=\alpha_1^{\beta_1}-\alpha_2^{\beta_2}$. Then, $\alpha_1$ and $\alpha_2$ satisfy $(\ref{2.14})$. By $(\ref{1.1})$ and $(\ref{2.23})$, we have $\Lambda^\prime=(-1)^{(a-1)x/2}(a^x+b^y)=(-1)^{(a-1)x/2}c^z$. It implies that $\Lambda^\prime\not=0$ and

 \begin{equation}\label{2.24}
{\rm ord_2}\Lambda^\prime\geq z.
\end{equation}
Further, using Lemma $\ref{2l3}$, by $(\ref{2.24})$, we have
 $$z<19.57(\log a)(\log b)(\max\{12\log2,$$
\begin{equation}\label{2.25}
0.4+\log(2\log2)+\log(\frac{x}{\log b}+\frac{y}{\log a})\})^2.
\end{equation}

 When $12\log2\geq0.4+\log(2\log2)+\log(x/\log b+y/\log a)$, by  $(\ref{2.25})$, we have

 \begin{equation}\label{2.26}
z<1355(\log a)(\log b).
\end{equation}
Further, by $(\ref{2.2})$ and $(\ref{2.26})$, we get $x<1355(\log b)(\log c)$ and $y<1355(\log a)(\log c)$, hence, $(\ref{2.1})$ holds.

When $12\log2<0.4+\log(2\log2)+\log(x/\log b+y/\log a)$, we get

\begin{equation}\label{2.27}
z<19.57(\log a)(\log b)(0.74+\log(\frac{x}{\log b}+\frac{y}{\log a}))^2.
\end{equation}
By $(\ref{2.2})$, we have $\max\{x/\log b,y/\log a\}<z(\log c)/(\log a)(\log b)$. Hence, by $(\ref{2.27})$, we get
\begin{equation}\label{2.28}
\frac{z\log c}{(\log a)(\log b)}<19.57(\log c)(1.44+\log(\frac{z\log c}{(\log a)(\log b)}))^2.
\end{equation}
Let $t=z(\log c)/(\log a)(\log b)$ and $F(t)=t-19.57(\log c)(1.44+\log t)^2$. Then,we have $F(3000(\log c)^2)>0$ and $F^\prime(t)=1-39.14(\log c)(1.44+\log t)/t>0$ for $t\geq3000(\log c)^2$. Hence, we have $F(t)>0$ for $t\geq 3000(\log c)^2$. Therefore, we see from $(\ref{2.28})$ that
 \begin{equation}\label{2.29}
z<3000(\log a)(\log b)(\log c).
\end{equation}
Further, by $(\ref{2.2})$ and $(\ref{2.29})$, we get $x<3000(\log b)(\log c)^2$ and $y<3000(\log a)(\log c)^2$, hence, $(\ref{2.1})$ holds. Thus, the theorem is true if $2|c$. To sum up, the theorem is proved.$\hfill{} \Box$
 \section {Proof of Theorem 1.1}

Let $m$ be a positive integer with $m>1$, and let $r,s$ be nonzero integers. The following lemma contains certain plain facts in elementary number theory.

 \begin{lemma}\label{3l1}.
\begin{enumerate}
\rm \item $\min\{|r|,|s|\}\geq\gcd(r,s)$.
\rm \item {\it If $r\equiv s\pmod m$, then $\gcd(r,m)=\gcd(s,m)$}.
\rm \item {\it If $\gcd(r,m)=1$, then $\gcd(rs,m)=\gcd(s,m)$}.
\rm \item {\it If $r|s$, then $\gcd(r,m)|\gcd(s,m)$}.
\rm \item {\it If $r>0,s>0$ and $rs\equiv0\pmod m$, then $s\geq m/\gcd(r,m)$}.
\rm \item {\it If $\gcd(r,m)=1$, then there exist integers $\bar{r}$ such that $r\bar{r}\equiv1\pmod m$ and $\gcd(\bar{r},m)=1$.}
\end{enumerate}
\end{lemma}

 \begin{lemma}\label{3l2}{\rm (\cite{car4})}.
If $\gcd(r,m)=1$, then there exist positive integers $n$ such that
\begin{equation}\label{3.1}
r^n\equiv\delta \pmod m,\delta\in\{1,-1\}.
\end{equation}
Let $n_1$ be the least value of $n$ with $(\ref{3.1})$, and let $r^{n_1}\equiv \delta_1\pmod m$, where $\delta_1\in\{1,-1\}$. A positive integer $n$ satisfies $(\ref{3.1})$ if and only if $n_1|n$. Moreover, if $n_1|n$ and $r^{n_1}-\delta_1\not=0$, then $r^{n_1}-\delta_1|r^n-\delta$.
\end{lemma}

Let $A,B,k$ be fixed positive integers such that $\min\{A,B,k\}>1$ and $\gcd(A,B)=1$.

\begin{lemma}\label{3l3}{\rm (\cite{ben1})}.
The equation
\begin{equation}\label{3.2}
A^m-B^n=k,m,n\in\mathbb{N}
\end{equation}
has at most two solutions $(m,n)$.
\end{lemma}

\begin{lemma}\label{3l4}.
The equation
\begin{equation}\label{3.3}
A^m+B^n=k,m,n\in\mathbb{N}
\end{equation}
has at most two solutions $(m,n)$.
\end{lemma}

{\bf {\it Proof .}}\,\,We now assume that $(\ref{3.3})$ has three solutions $(m_i,n_i)(i=1,2,3)$. By the symmetry of $A^m$ and $B^n$ in $(\ref{3.3})$, we may therefore assume that $m_1>m_2$ and $A^{m_j}>B^{n_j}$ for $j=1,2$. Then, since $A^{m_1}+B^{n_1}=A^{m_2}+B^{n_2}$, we have $n_1<n_2$ and $A^{m_j}\equiv-B^{n_j}\pmod k$ for$j=1,2$. Hence, we get

\begin{equation}\label{3.4}
A^{m_1}B^{n_2}\equiv A^{m_2}B^{n_1}\pmod k.
\end{equation}
Further, since $\gcd(A,B)=1$, we have $\gcd(AB,k)=1$, and by $(\ref{3.4})$,

\begin{equation}\label{3.5}
A^{m_1-m_2}B^{n_2-n_1}\equiv 1\pmod k.
\end{equation}
It implies that $A^{m_1-m_2}B^{n_2-n_1}\geq k+1>k>A^{m_1}$, whence we obtain $B^{n_2}>B^{n_2-n_1}>A^{m_1-(m_1-m_2)}=A^{m_2}>B^{n_2}$, a contradiction. Thus, the lemma is proved.$\hfill{} \Box$\\


For any fixed triple $(a,b,c)$, put $P(a,b,c)=\{(a,b,c,1),(c,a,b,-1),(c,b,a,-1)\}$. Obviously, there exists a unique element in $P(a,b,c)$, say $(A,B,C,\lambda)$, which satisfies

\begin{equation}\label{3.6}
C=\max\{a,b,c\}.
\end{equation}
Then, $(\ref{1.1})$ has a solution $(x,y,z)$ is equivalent to the equation

\begin{equation}\label{3.7}
A^X+\lambda B^Y=C^Z, X,Y,Z\in\mathbb{N}
\end{equation}
has the solution

$$(X,Y,Z)=\left\{\begin{array}{cc}
(x,y,z),  &{{\rm if}\ \ (A,B,C,\lambda)=(a,b,c,1) ,} \\
(z,x,y),  &{{\rm if}\ \ (A,B,C,\lambda)=(c,a,b,-1),}\\
(z,y,x),&{{\rm if}\ \ (A,B,C,\lambda)=(c,b,a,-1).}
\end{array} \right.$$
It implies that the numbers of solutions of $(\ref{1.1})$ and $(\ref{3.7})$ are equal.

Here and below, we always assume that $(\ref{3.7})$ has solutions $(X,Y,Z)$. Then, it has a solution $(X_1,Y_1,Z_1)$ such that $Z_1\leq Z$, where $Z$ through all solutions $(X,Y,Z)$ of $(\ref{3.7})$. Since $\gcd(A,C)=1$, by Lemma $\ref{3l2}$, there exist positive integers $n$ such that

\begin{equation}\label{3.8}
A^n\equiv \delta\pmod{C^{Z_1}},\delta\in\{1,-1\}.
\end{equation}
Let $n_1$ be the least value of $n$ with $(\ref{3.8})$, and let
\begin{equation}\label{3.9}
A^{n_1}\equiv \delta_1\pmod{C^{Z_1}},\delta_1\in\{1,-1\}.
\end{equation}
Then we have
\begin{equation}\label{3.10}
A^{n_1}=C^{Z_1}f+\delta_1,f\in\mathbb{N}.
\end{equation}
Obviously, for any fixed triple $(a,b,c)$, the parameters $Z_1,n_1,\delta_1$ and $f$ are unique.

\begin{lemma}\label{3l5}{\rm (\cite{hu8}, Lemma 3.3)}.
If $(X,Y,Z)$ and $(X^\prime,Y^\prime,Z^\prime)$ are two solutions of $(\ref{3.7})$ with $Z\leq Z^\prime$, then $XY^\prime-X^\prime Y\not=0$ and
$$A^{|XY^\prime-X^\prime Y|}\equiv(-\lambda)^{Y+Y^\prime}\pmod {C^Z}.$$
\end{lemma}

\begin{lemma}\label{3l6}.
$(\ref{3.7})$ has at most two solutions $(X,Y,Z)$  with $Z=Z_1$ .
\end{lemma}
{\bf {\it Proof .}}\,\, By Lemmas $\ref{3l3}$ and $\ref{3l4}$, we obtain the lemma immediately.$\hfill{} \Box$\\

\begin{lemma}\label{3l7}.
If $(\ref{3.7})$ has two solutions $(X_1,Y_1,Z_1)$  and $(X_2,Y_2,Z_2)$ with $Z_1<Z_2$, then $\gcd(C,f)\leq Y_2$.
\end{lemma}
{\bf {\it Proof .}}\,\, Since $A^{X_1}+\lambda B^{Y_1}=C^{Z_1}, A^{X_2}+\lambda B^{Y_2}=C^{Z_2}$ and $Z_1+1\leq Z_2$, we have
$$A^{X_1Y_2}\equiv(-\lambda)^{Y_2}B^{Y_1Y_2}+(-\lambda)^{Y_2-1}B^{Y_1(Y_2-1)}C^{Z_1}Y_2\pmod {C^{Z_1+1}},$$
\begin{equation}\label{3.11}
A^{X_2Y_1}\equiv(-\lambda B^{Y_2}+C^{Z_2})^{Y_1}\equiv(-\lambda)^{Y_1}B^{Y_1Y_2}\pmod {C^{Z_1+1}}.
\end{equation}
Eliminating $B^{Y_1Y_2}$ from $(\ref{3.11})$, we get
$$A^{\min\{X_1Y_2,X_2Y_1\}}(A^{|X_1Y_2-X_2Y_1|}-(-\lambda )^{Y_1+Y_2})$$
\begin{equation}\label{3.12}
\equiv \lambda^\prime B^{Y_1(Y_2-1)}C^{Z_1}Y_2\pmod {C^{Z_1+1}},
\end{equation}
where
$$\lambda^\prime=\left\{\begin{array}{cc}
(-\lambda)^{Y_2-1},  &{{\rm if}\ \ X_1Y_2>X_2Y_1,} \\
-(-\lambda)^{Y_1-1},  &{{\rm if}\ \ X_1Y_2<X_2Y_1.}
\end{array} \right.$$
Further, since $\gcd(A,C)=1$, by (\romannumeral6) of Lemma $\ref{3l1}$, there exist integers $\bar{A}$ such that $A\bar{A}\equiv1\pmod {C^{Z_1+1}}$ and $\gcd(\bar{A},C)=1$. Hence, by $(\ref{3.12})$, we have
$$A^{|X_1Y_2-X_2Y_1|}-(-\lambda )^{Y_1+Y_2}$$
\begin{equation}\label{3.13}
\equiv\lambda^\prime\bar{A} ^{\min\{X_1Y_2,X_2Y_1\}}B^{Y_1(Y_2-1)}C^{Z_1}Y_2\pmod {C^{Z_1+1}}.
\end{equation}

By Lemma $\ref{3l5}$, $|X_1Y_2-X_2Y_1|$ is a positive integer. We see from $(\ref{3.13})$ that

\begin{equation}\label{3.14}
A^{|X_1Y_2-X_2Y_1|}-(-\lambda )^{Y_1+Y_2}= C^{Z_1}g, g\in\mathbb {N},
\end{equation}
where $g$ satisfies

\begin{equation}\label{3.15}
g\equiv \lambda ^\prime \bar{A} ^{\min\{X_1Y_2,X_2Y_1\}}B^{Y_1(Y_2-1)}Y_2\pmod C.
\end{equation}
Applying (\romannumeral2) of Lemma $\ref{3l1}$ to $(\ref{3.15})$, we get

\begin{equation}\label{3.16}
\gcd(C,g)=\gcd(C,\lambda ^\prime \bar{A} ^{\min\{X_1Y_2,X_2Y_1\}}B^{Y_1(Y_2-1)}Y_2).
\end{equation}
Further, since $\lambda^\prime\in\{1,-1\}$ and $\gcd(\bar{A},C)=\gcd(B,C)=1$, by
(\romannumeral3) of Lemma $\ref{3l1}$, we have
\begin{equation}\label{3.17}
\gcd(C,\lambda ^\prime \bar{A} ^{\min\{X_1Y_2,X_2Y_1\}}B^{Y_1(Y_2-1)}Y_2)=\gcd(C,Y_2).
\end{equation}
The combination of $(\ref{3.16})$ and $(\ref{3.17})$ yields
\begin{equation}\label{3.18}
\gcd(C,g)=\gcd(C,Y_2).
\end{equation}

On the other hand, by $(\ref{3.14})$, we have
\begin{equation}\label{3.19}
A^{|X_1Y_2-X_2Y_1|}\equiv(-\lambda )^{Y_1+Y_2}\pmod {C^{Z_1}}.
\end{equation}
Applying Lemma $\ref{3l2}$ to  $(\ref{3.9})$ and $(\ref{3.19})$, we get $n_1|X_1Y_2-X_2Y_1$ and $A^{n_1}-\delta_1|A^{|X_1Y_2-X_2Y_1|}-(-\lambda)^{Y_1+Y_2}$. Hence, by  $(\ref{3.10})$ and $(\ref{3.14})$, we have
\begin{equation}\label{3.20}
f|g.
\end{equation}
Therefore, using (\romannumeral4) of Lemma $\ref{3l1}$, by $(\ref{3.20})$, we get $\gcd(C,f)|\gcd(C,g)$ and
\begin{equation}\label{3.21}
\gcd(C,f)\leq\gcd(C,g).
\end{equation}
Further, by (\romannumeral1) of Lemma $\ref{3l1}$, we have $\gcd(C,Y_2)\leq Y_2$. Thus, by $(\ref{3.18})$ and $(\ref{3.21})$, we obtain $\gcd(C,f)\leq Y_2$. The lemma is proved.$\hfill{} \Box$\\

\begin{lemma}\label{3l8}.
If $(\ref{3.7})$ has three solutions $(X_i,Y_i,Z_i)(i=1,2,3)$  with $Z_1<Z_2\leq Z_3$, then $\max\{a,b,c\}<5\times10^{27}$.
\end{lemma}
{\bf {\it Proof .}}\,\,Since $Z_1+1\leq Z_2\leq Z_3$, by Lemma $\ref{3l5}$, we have $X_2Y_3-X_3Y_2\not=0$ and

\begin{equation}\label{3.22}
A^{|X_2Y_3-X_3Y_2|}\equiv(-\lambda )^{Y_2+Y_3}\pmod {C^{Z_1+1}}.
\end{equation}
It implies that
\begin{equation}\label{3.23}
A^{|X_2Y_3-X_3Y_2|}-(-\lambda )^{Y_2+Y_3}=C^{Z_1+1}h, h\in\mathbb{N}.
\end{equation}

On the other hand, using Lemma $\ref{3l2}$, we see from $(\ref{3.9})$ and $(\ref{3.22})$ that $n_1|X_2Y_3-X_3Y_2$ and

\begin{equation}\label{3.24}
|X_2Y_3-X_3Y_2|=n_1n_2, n_2\in\mathbb{N}.
\end{equation}
By $(\ref{3.10})$,  $(\ref{3.23})$ and $(\ref{3.24})$, we have

$$C^{Z_1+1}h=A^{n_1n_2}-(-\lambda)^{Y_2+Y_3}=(C^{Z_1}f+\delta_1)^{n_2}-(-\lambda)^{Y_2+Y_3}$$

\begin{equation}\label{3.25}
=(\delta_1^{n_2}-(-\lambda)^{Y_2+Y_3})+C^{Z_1}f\sum\limits_{i=1}^{n_2}\left(\begin{array}{cc}
n_2 \\
i
\end{array} \right)\delta_1^{n_2-i}(C^{Z_1}f)^{i-1}.
\end{equation}
Since $C>2$ by $(\ref{3.6})$, we find from $(\ref{3.25})$ that $\delta_1^{n_2}=(-\lambda)^{Y_2+Y_3}$ and

\begin{equation}\label{3.26}
Ch=f\sum\limits_{i=1}^{n_2}\left(\begin{array}{cc}
n_2 \\
i
\end{array} \right)\delta_1^{n_2-i}(C^{Z_1}f)^{i-1},
\end{equation}
whence we get
\begin{equation}\label{3.27}
fn_2\equiv0\pmod C.
\end{equation}

Applying (\romannumeral5) of Lemma $\ref{3l1}$ to $(\ref{3.27})$, we have

\begin{equation}\label{3.28}
n_2\gcd(C,f)\geq C.
\end{equation}
Further, by  Lemma $\ref{3l7}$ and $(\ref{3.24})$, we have $\gcd(C,f)\leq Y_2$ and $n_2\leq|X_2Y_3-X_3Y_2|<\max\{X_2Y_3,X_3Y_2\}$ respectively. Therefore, by $(\ref{3.28})$, we get

\begin{equation}\label{3.29}
C<Y_2\max\{X_2Y_3,X_3Y_2\}\leq(\max\{X_2,Y_3,X_3,Y_2\})^3.
\end{equation}
Recall that every solution $(X,Y,Z)$ of $(\ref{3.7})$ is a permutation of a solution $(x,y,z)$ of $(\ref{1.1})$. By Theorem 2.1, we have $$\max\{X_2,Y_2,X_3,Y_3\}<6500(\log\max\{a,b,c\})^3.$$
Hence, by $(\ref{3.6})$ and $(\ref{3.29})$, we get

\begin{equation}\label{3.30}
\max\{a,b,c\}<6500^3(\log\max\{a,b,c\})^9.
\end{equation}

Let $t=\max\{a,b,c\}$ and $F(t)=t-6500^3(\log t)^9$. Then we have $F(5\times10^{27})>0$ and $F^\prime(t)=1-9\times6500^3(\log t)^8/t>0$ for $t\geq5\times10^{27}$. Therefore, we get $F(t)>0$ for $t\geq 5\times10^{27}$. Thus, we obtain from $(\ref{3.30})$ that $\max\{a,b,c\}<5\times10^{27}$. The lemma is proved.$\hfill{} \Box$\\

{\bf {\it Proof of Theorem 1.1.}}\,\,We now assume that $(\ref{1.1})$ has four solutions $(x,y,z)$. Then $(\ref{3.7})$ has four solutions $(X,Y,Z)$. Further, by Lemma $\ref{3l6}$, $(\ref{3.7})$ has three solutions $(X_i,Y_i,Z_i)(i=1,2,3)$ with $Z_1<Z_2\leq Z_3$. Therefore, by Lemma $\ref{3l8}$, we get $\max\{a,b,c\}<5\times10^{27}$. Thus, the theorem is proved.$\hfill{} \Box$

\begin{flushleft}
Yongzhong Hu\\
Department of Mathematics \\
Foshan University\\
Foshan,Guangdong 528000,China\\
E-mail:huuyz@aliyun.com
\end{flushleft}

\begin{flushleft}
Maohua Le\\
Institute of Mathematics\\
Lingnan Normal University\\
Zhanjiang,Guangdong 524048,China\\
E-mail:lemaohua2008@163.com
\end{flushleft}

    \end{document}